   \documentclass[11pt]{amsart}
\usepackage{amsmath,amssymb,amscd,amsthm,eucal,curves,epic, graphicx}
\usepackage[all]{xy}
\usepackage[latin1]{inputenc}
\usepackage{color}
\usepackage{floatflt}
 
\newcommand{\Q}{{\mathbb Q}}

\newcommand{\Map}{\operatorname{Map}_*}

\newcommand{\colim}{\operatorname{colim}}

\newcommand{\Spaces}{\bf{Spaces}_*}

\newtheorem{theorem}{Theorem}[section]

\newtheorem{corollary}[theorem]{Corollary}

\theoremstyle{definition}

\newtheorem{definition}[theorem]{Definition}

\input xy
\xyoption{all}

\title[Localizations of mapping spaces]{A note on localizations of mapping spaces}
\date{07/15/2008}                                       

\author[B. Badzioch]{Bernard Badzioch}
\author[W. Dorabia{\l}a]{Wojciech Dorabia{\l}a}

\address[]{Department of Mathematics, University at Buffalo, SUNY, Buffalo, NY}
\address[]{Department of Mathematics, Penn State Altoona, Altoona, PA }

\begin{document}

\begin{abstract}
We show that if $A$ is a simply connected, finite, pointed CW-complex then the
mapping spaces $\Map(A, X)$ are preserved by the localization functors only 
if $A$ has the rational homotopy type of a wedge of spheres $\bigvee_{l}S^{k}$.
\end{abstract}

\maketitle

%
%
 
 \section{Introduction}
 \label{INTRO}
 The motivation for this brief note comes from the following well known property of localization 
 functors \cite[Thm 3.A.2]{FARJOUN}. Given a map of pointed spaces $f$ consider the localization
 functor $L_{f}\colon \Spaces\to \Spaces$.  For any $X\in \Spaces$  we have a weak equivalence
 \begin{equation}
 \label{Lf EQ}
 L_{f}\Omega X\simeq \Omega L_{\Sigma f}X
 \end{equation}
 This shows that localizations preserve loop spaces. 
  
It is natural to ask if this preservation property can be extended. This leads to the following
\begin{definition}
\label{L-GOOD}
We say that a finite,  connected, pointed CW-complex $A$ is $L$-good if for any pointed map 
$f$ and any $X\in\Spaces$ we have
$$L_{f}\Map(A, X)\simeq \Map(A, Y)$$
for some $Y\in \Spaces$.
\end{definition}
The weak equivalence (\ref{Lf EQ}) shows that $S^{1}$ is $L$-good. We would like to know what 
other spaces have this property.  This is in fact one of the questions posed by Dror Farjoun in 
\cite[9.F]{FARJOUN}.
 Since $\Omega^{k}X\cong \Omega(\Omega^{k-1}X)$, applying iteratively the weak equivalence 
 (\ref{Lf EQ}) we get that $S^{k}$ is $L$-good 
for all $k\geq 1$. Also, since $\Map(\bigvee_{l}S^{k}, X)\cong \prod^{l}\Map(S^{k}, X)$,
and since localization functors preserve finite products up to a weak equivalence, 
we obtain that the class of $L$-good spaces contains all spaces $\bigvee_{l}S^{k}$ for $k>0$, 
$l\geq 0$. Our goal here is to show that, rationally, every $L$-good space will resemble 
$\bigvee_{l}S^{k}$. 

\begin{theorem}
\label{RATIONAL THM}
Let $A$ be a finite, connected, pointed CW-complex such that for some $p>q>0$ we have 
$H^{p}(A, \Q)\neq 0\neq H^{q}(A, \Q)$. Then $A$ is not an $L$-good space. 
 \end{theorem}

Equivalently,  for an $L$-good space $A$ we have $H^{i}(A, \Q)\neq 0$ for at most one
 $i>0$. As a consequence we obtain
 
 \begin{corollary}
 If $A$ is a simply connected $L$-good space then $A$ has the rational homotopy type 
 of $\bigvee_{l}S^{k}$ for some $k>0$, $l\geq0$.
 \end{corollary}
 
 We note here that the formula (\ref{Lf EQ}) follows from the existence of 
 the loop space machines (see e.g.  \cite{BOARDMAN},  \cite{MAY}, \cite{SEGAL}) 
 which describe the structure of spaces $\Omega X$ in 
 terms of maps of finite products $(\Omega X)^{m} \to (\Omega X)^{n}$. An analogous 
 description of mapping spaces $\Map(A, X)$ for some $A$ would simi\-larly imply that $A$ is an
 $L$-good space. Theorem \ref{RATIONAL THM} shows then that finite product 
 ``mapping space'' machines do not exist for any finite CW-complex $A$ whose rational 
 cohomology is non-trivial in more than one dimension.
 
 {\bf Acknowledgements.} This work was completed while the authors par\-ti\-cipa\-ted in 
 the Research in Pairs program at the Mathematisches For\-schungs\-institut Oberwolfach. 
 The authors want to express their gratitute to the Institut for its hospitality.
 The first author also wants to thank A. Prze\'zdziecki and W. G. Dwyer for conversations which 
 inspired this paper.

%
%

 \section{Proof of Theorem \ref{RATIONAL THM}}
 \label{RATIONAL PROOF}
 
Let $A$ be a CW-complex as in the statement of Theorem \ref{RATIONAL THM}. 
Since $A$ is finite we can choose $p$ so that $H^{i}(A, \Q)=0$ for all $i>p$.
For $n>p$ we have a weak equivalence
$$\Map(A, K(\Q, n))\simeq \prod_{i=n-p}^{n}K(H^{n-i}(A, \Q), i)$$
Consider the constant map  $f\colon S^{k}\to \ast$. In this case the localization 
$L_{f}$ is the nullification functor $P_{S^{k}}$. We have
$$P_{S^{n-p+1}}\Map(A, K(\Q, n))\simeq K(H^{p}(A, \Q), n-p)$$
If follows that if $A$ was an $L$-good space then for every $N>0$ we would 
be able to find a space $Y$ such that 
\begin{equation}
\label{Y EQ}
\Map(A, Y)\simeq K(H^{p}(A, \Q), N)
\end{equation}

We will show that this is impossible arguing by contradiction.  
Assume first that $A$ is simply connected, $0\neq V =H^{p}(A, \Q)$, and that for some fixed 
$N>p+1$ we have a space $Y$ satisfying (\ref{Y EQ}).

Since $A$ is simply connected we have 
$\Map(A, Y)\simeq \Map(A, \widetilde{Y})$ where $\widetilde{Y}$ is the universal cover
of $Y$. Therefore we can assume that $Y$ is simply connected. 

Next, let  $Y_{(0)}$ denote the rationalization of $Y$. By \cite[Thm.3.11, p.77]{HILTON}
$\Map(A, Y_{(0)})\simeq \Map(A, Y)_{(0)}$, and since $\Map(A, Y)\simeq K(V, N)$ 
is a rational space thus $\Map(A, Y_{(0)})\simeq \Map(A, Y)$.  As a consequence 
we can assume that $Y$ is a simply connected rational space. 

By \cite[Corollary p. 229]{FELIX} we have 
$$\Omega Y \simeq \tilde \prod_{n\geq 1}K(V_{n}, n)$$
where $V_{n}$ is a $\Q$-vector space and $\tilde\prod$ denotes the weak product of pointed 
spaces: $\tilde \prod_{n\geq 1} K(V_{n}, n)=\colim_{M\geq 1} (\prod_{n=1}^{M}K(V_{n}, n))$.
We obtain

\begin{eqnarray}
\label{EQ N-1}
K(V, N-1)\simeq \Map(A,\Omega Y)\simeq 
\Map(A,  \tilde \prod_{n\geq 1}K(V_{n}, n))
\end{eqnarray}

We claim that there exists $n_{0}\geq N-1$ such that $V_{n_{0}}\neq 0$.  Indeed, if
$V_{n}=0$ for all $n\geq N-1$ then 
$ \tilde \prod_{n\geq 1}K(V_{n}, n)=\prod_{n=1}^{N-2}K(V_{n}, n)$ so
 
$$\Map(A,  \tilde \prod_{n\geq 1}K(V_{n}, n))=\prod_{n=1}^{N-2}\Map(A, K(V_{n}, n))$$
This would give 
$$\pi_{i}(\Map(A,  \tilde \prod_{n\geq 1}K(V_{n}, n)))\cong \bigoplus_{n=1}^{N-2}
\widetilde H^{n-i}(A, V_{n})$$ 
In particular we would have $\pi_{i}(\Map(A,  \tilde \prod_{n\geq 1}K(V_{n}, n)))=0$ for $i\geq N-1$
which contradicts  (\ref{EQ N-1}).

Since $n_{0}\geq N-1> p, q$ we have 
$$\pi_{n_{0}-p}(\Map(A, K(V_{n_{0}}, n_{0}))\cong H^{p}(A, V_{n_{0}})\neq 0$$ 
and
$$\pi_{n_{0}-q}(\Map(A, K(V_{n_{0}}, n_{0}))\cong  H^{q}(A, V_{n_{0}})\neq 0$$
where the inequalities on the right hold by our assumption that   
$H^{p}(A, \Q)\neq 0$, $H^{q}(A, \Q)\neq 0$. Also, the space $\Map(A, K(V_{n_{0}}, n_{0}))$ is a retract 
of $\Map(A,  \tilde \prod_{n\geq 1}K(V_{n}, n))$ so this last space must have non-trivial 
homotopy groups in at least two dimensions $n_{0}-p$ and $n_{0}-q$. This however contradicts
the formula (\ref{EQ N-1}). The contradiction shows that $\Map(A, Y)\not\simeq K(V, N)$ for any 
space $Y$, and so $A$ is not an $L$-good space. 

Assume now that $A$ is not simply connected. If $A$ was an $L$-good space then again we would 
be able to find a space $Y$ such that $\Map(A, Y)\simeq K(V, N)$, where $V=H^{p}(A, \Q)$, 
$N>p+2$. This would give
$$\Map(\Sigma A, Y)\simeq \Omega\Map(A, Y)\simeq K(V, N-1)$$
Since $\Sigma A$ is a simply connected space this is however impossible by the argument above. 
It follows that $\Map(A, Y)\not \simeq K(V, N)$ for any $Y\in \Spaces$, and so $A$ is not an 
$L$-good space.

\bibliographystyle{plain}
\bibliography{localizations}

\begin{thebibliography}{1}

\bibitem{BOARDMAN}
J.~M. Boardman and R.~M. Vogt.
\newblock {\em Homotopy invariant algebraic structures on topological spaces}.
\newblock Lecture Notes in Mathematics, Vol. 347. Springer-Verlag, Berlin,
  1973.

\bibitem{FARJOUN}
Emmanuel~Dror Farjoun.
\newblock {\em Cellular spaces, null spaces and homotopy localization}, volume
  1622 of {\em Lecture Notes in Mathematics}.
\newblock Springer-Verlag, Berlin, 1996.

\bibitem{FELIX}
Yves F{\'e}lix, Stephen Halperin, and Jean-Claude Thomas.
\newblock {\em Rational homotopy theory}, volume 205 of {\em Graduate Texts in
  Mathematics}.
\newblock Springer-Verlag, New York, 2001.

\bibitem{HILTON}
Peter Hilton, Guido Mislin, and Joe Roitberg.
\newblock {\em Localization of nilpotent groups and spaces}.
\newblock North-Holland Publishing Co., Amsterdam, 1975.
\newblock North-Holland Mathematics Studies, No. 15, Notas de Matem\'atica, No.
  55. [Notes on Mathematics, No. 55].

\bibitem{MAY}
J.~P. May.
\newblock {\em The geometry of iterated loop spaces}.
\newblock Springer-Verlag, Berlin, 1972.
\newblock Lectures Notes in Mathematics, Vol. 271.

\bibitem{SEGAL}
Graeme Segal.
\newblock Categories and cohomology theories.
\newblock {\em Topology}, 13:293--312, 1974.

\end{thebibliography}

\end{document}